\documentclass[a4, 11pt]{amsart}
\usepackage{amsmath,amsthm,amsfonts,amssymb,amscd}

\newtheorem{theorem}{Theorem}[section]

\usepackage{xcolor}
\usepackage{mathrsfs}
\usepackage{blkarray}
\usepackage{centernot,tikz}

\newtheorem{definition}[theorem]{Definition}
\newtheorem{example}[theorem]{Example}

\newtheorem{prop}[theorem]{Proposition}
\newtheorem{corollary}[theorem]{Corollary}
\newtheorem{question}[theorem]{Question}
\newtheorem{remark}[theorem]{Remark}

\numberwithin{equation}{section}
\DeclareMathOperator{\spn}{span}
\begin{document}
\title[Conditions implying the normality of $\ast$-paranormal operators in $\overline{\mathcal{AN}(H)}$]{Conditions implying the normality of $\ast$-paranormal operators in the closure of $\mathcal{AN}$-operators}
\author{G. Ramesh}
\address{G. Ramesh, Department of Mathematics, IIT Hyderabad, Kandi, Sangareddy, Telangana- 502284, India.}
 \email{rameshg@math.iith.ac.in}
\author{Shanola S. Sequeira}
\address{Shanola S. Sequeira, Department of Mathematics, IIT Hyderabad, Kandi, Sangareddy, Telangana- 502284, India.}
\email{ma18resch11001@gmail.com}
\maketitle
\begin{abstract}
	In this article, we first prove the existence of an invariant subspace for a norm attaining $\ast$-paranormal operator. Then give a representation for $\ast$-paranormal operators in the closure of absolutely norm attaining operators and further study a few sufficient conditions for the normality of such operators. Finally, we discuss Toeplitz and Hankel $\ast$-paranormal operators in the closure of absolutely norm attaining operators on the Hardy space.
\end{abstract}

\section{Introduction}
One of the attractive areas of research in operator theory is the study of non-normal operators. The class of non-normal operators includes hyponormal, paranormal operators, etc., and are studied extensively in the literature \cite{BER1,Furutapara,Istratescu,STA}. Another important class of non-normal operators is the class of $\ast$-paranormal operators.  A bounded linear operator $T$ on a Hilbert space $H$ is said to be $\ast$-paranormal if
\[\|T^*x\|^2 \leq \|T^2x\| \|x\|, \ x \in  H. \]  S.C. Arora and J.K. Thukral initiated the study of $\ast$-paranormal operators \cite{Arora}. Later on, this class received a lot of attention with the development of Weyl's theorem \cite{HAN,Kim}. A study on more generalized classes containing $\ast$-paranormal operators is also done (see \cite{Rashidn*,RashidT*para,Uchi*para} for more details). In \cite{Uchi*para}, the authors gave an example of an invertible $\ast$-paranormal operator whose inverse is not
$\ast$-paranormal. But this is not true for paranormal operators, as the inverse of an invertible paranormal operator is paranormal \cite{Istratescu}. This shows that even though both classes contain hyponormal operators they do not coincide. This makes the study on $\ast$-paranormal operators as interesting as the study on other classes of non-normal operators.

It is well known that a compact hyponormal operator is normal \cite{BER1,STA}, as well as a compact paranormal operator is normal \cite{QIU}. Further, a compact $\ast$-paranormal operator is also normal \cite{Rashidn*}. That means under the assumption of compactness, these non-normal operators become normal. Hence it is natural to ask the following question:
\begin{question}\label{normalityquestion}
	Under what assumptions, does a non-normal operator become normal?
\end{question}  This question is also motivated by the results of Berberian \cite{BER1}, Putnam \cite{Putnam}, Qiu \cite{QIU}, and Stampfli \cite{STA}, where they discussed some conditions implying the normality of paranormal and hyponormal operators, which include the countability of the spectrum, spectrum with area measure zero, etc. Hence it is interesting to answer Question \ref{normalityquestion}
for $\ast$-paranormal operators also by replacing compactness with a weaker property namely absolutely norm attaining property.

Throughout the article, let $H, H_1, H_2$ denote infinite dimensional complex Hilbert spaces and $\mathcal{B}(H_1, H_2)$ denote the Banach space of all bounded linear operators from $H_1$ into $H_2$. If $H_1= H_2 =H$, then we denote $\mathcal{B}(H_1, H_2)$ by $\mathcal{B}(H)$.\begin{definition}\cite[Definitions 1.1, 1.2]{CAR1}
	An operator $T \in \mathcal{B}(H_1,H_2)$ is called norm attaining
	if there exists a unit vector $x \in H_1$ such that $\|T\| = \|Tx\|$. If $T|_M :M \to H_2$  is norm attaining for every non-zero closed subspace $M$ of $H_1$,  then $T$ is called an absolutely norm attaining or $\mathcal{AN}$-operator.
\end{definition}
This class contains compact operators, isometries and partial isometries with finite dimensional null space.
The set of norm attaining and absolutely norm attaining operators from $H_1$ to $H_2$ are denoted by $\mathcal{N}(H_1, H_2)$ and $\mathcal{AN}(H_1, H_2)$, respectively. We denote $\mathcal{N}(H, H)$ by $\mathcal{N}(H)$ and $\mathcal{AN}(H, H)$ by $\mathcal{AN}(H)$. For a recent account of the theory of $\mathcal{AN}$-operators, we refer to \cite{CAR1,PAN,Rameshpara,VEN}. The right shift operator on $l^2$ is $\ast$-paranormal and an isometry. So every absolutely norm attaining $\ast$-paranormal operator need not be normal. Hence the author in \cite{NB*-para} was successful in giving a representation of $\ast$-paranormal $\mathcal{AN}$-operators, and some conditions under which they become normal.

Similar to the concept of the norm of an operator, we have the concept of minimum modulus which is defined by
\begin{equation*}
m(T):= \inf\{\|Tx\| : x \in H_1, \|x\| = 1 \}.
\end{equation*}
Analogous to norm attaining and absolutely norm attaining operators, we have the following classes of operators.
\begin{definition}\cite[Definitions 1.1, 1.4]{CAR2}
	An operator $T \in \mathcal{B}(H_1, H_2)$ is called minimum attaining if there exists a unit vector $x \in H_1$ such that $m(T) = \|Tx\|$. If $T|_M: M \to H_2$  is minimum attaining for every non-zero closed subspace $M$ of $H_1$,  then $T$ is called an absolutely minimum attaining or $\mathcal{AM}$-operator.
\end{definition}
The sets of all minimum attainining and absolutely minimum attatining operators from $H_1$ to $H_2$ are denoted by $\mathcal{M}(H_1,H_2)$ and $\mathcal{AM}(H_1,H_2)$, respectively and when $H_1=H_2=H$, we write $\mathcal{M}(H,H):= \mathcal{M}(H)$ and $\mathcal{AM}(H,H):= \mathcal{AM}(H)$. This class contains finite rank operators, isometries and partial isometries with finite dimensional null space. We refer to \cite{BALA1,CAR2,GAN} for more information on these operators.


In this article, we examine Question \ref{normalityquestion} with respect to a larger class, namely the operator norm closure of $\mathcal{AN}$-operators. This class includes both $\mathcal{AN}$-operators and $\mathcal{AM}$-operators. In addition, $\mathcal{AN}$-operators have the same operator norm closure as $\mathcal{AM}$-operators. Also, the closure of $\mathcal{AN}$-operators neither contains nor is contained by the class of norm attaining operators. All of these details are available in \cite{RAMSSS}. In addition, this class contains operators whose essential spectrum is a singleton set, such as quasinilpotent operators and Gaussian Covariance operators \cite{Bhat}. The results of this article generalise the conclusions of \cite{NB*-para}, however the methodologies employed in\cite{NB*-para} are different.


In the first step of this investigation, we demonstrate the existence of an invariant subspace for a norm attaining $\ast$-paranormal operator distinct from the one discussed in \cite{NB*-para}. Then, we offer a representation of a $\ast$-paranormal operator in $\overline{\mathcal{AN}(H)}$ that generalises the representation of hyponormal operators in \cite{RAMSSS1}. As a result, we obtain the following.
\begin{enumerate}
\item Representation of  $\ast$-paranormal $\mathcal{AM}$-operators.
\item Representation of  $\ast$-paranormal $\mathcal{AM}$-operators.
\end{enumerate}

We further show that any invertible $\ast$-paranormal operator in $\overline{\mathcal{AN}(H)}$ is normal and give a few more conditions for normality. All these results are discussed in section 2.

Finally, in section 3, we show that a $\ast$-paranormal Toeplitz operator in $\overline{\mathcal{AN}(H^2)}$ is a scalar multiple of an isometry and a $\ast$-paranormal Hankel operator in $\overline{\mathcal{AN}(H^2)}$ is normal, where $H^2$ is the Hardy space of the unit circle of the complex plane.

In the remaining part of this section, we give all the necessary definitions and results required to develop the article.

 \subsection{Preliminaries}
For any $T \in \mathcal{B}(H_1, H_2)$, the adjoint operator $T^* \in \mathcal{B}(H_2, H_1)$ is defined by
\[\langle Tx, y \rangle = \langle x, T^{*}y \rangle, \ \forall x \in H_1,\  y \in H_2.\]
Now, we define the notion of the spectrum of an operator which generalizes the concept of eigenvalues of a matrix. If $T \in \mathcal{B}(H)$, then the spectrum is defined by
\[\sigma(T) \ := \{\lambda \in \mathbb{C} : T - \lambda I \ \text{is} \ \text{not} \ \text{invertible} \ \text{in} \ \mathcal{B}(H)\}.\]

The spectrum of $T$ is decomposed as the disjoint union of the point spectrum $\sigma_{p}(T):= \{\lambda \in \mathbb{C} : T-\lambda I \ \text{is}\; \text{not}\; \text{one-one} \ \text{in} \ \mathcal{B}(H)\}$, the residual spectrum $\sigma_{r}(T) := \{\lambda \in \mathbb{C} : T-\lambda I \ \text{is}\; \text{one-one}\; \text{but}\; \overline{R(T-\lambda I)} \neq H\}$ and the continuous spectrum $\sigma_{c}(T) := \sigma(T) \setminus (\sigma_{p}(T) \cup \sigma_{r}(T))$.

We say $T \in \mathcal{B}(H)$ is normal if $TT^* = T^*T$, self-adjoint if $T = T^*$, and positive if it is self-adjoint and $\langle Tx, x\rangle \geq 0$. A positive operator $T \in \mathcal{B}(H)$ is denoted by $T \geq 0$.

If $T \in \mathcal{B}(H_1, H_2)$, we call $|T| = (T^*T)^{1/2}$ as the modulus of $T$. For $T \in \mathcal{B}(H_1, H_2)$ there exists a unique partial isometry $W \in \mathcal{B}(H_1, H_2)$ with $N(W) = N(T)$ such that $T = W |T|$. This is called the polar decomposition of $T$.

Let $R(T)$ and $N(T)$ be the range and nullspaces of $T \in \mathcal{B}(H_1, H_2)$, respectively. We say $T$ to be a finite rank operator if $R(T)$ is finite dimensional and $T$ to be a compact operator if, for any bounded set $B$ of $H_1$, $T(B)$ has a compact closure. We denote the set of all finite rank and compact operators from $H_1$ to $H_2$ by $\mathcal{F}(H_1, H_2)$ and $\mathcal{K}(H_1, H_2)$, respectively.  If $H_1 = H_2 = H$, then $\mathcal{K}(H_1, H_2) = \mathcal{K}(H)$ and $\mathcal{F}(H_1, H_2) = \mathcal{F}(H)$.

If $M$ is a closed subspace of $H$, then $M^\bot$ denotes the orthogonal complement of $M$, $P_M$ denotes the orthogonal projection onto $M$ and $I_M$ denotes the identity operator on $M$, respectively. Let $T \in \mathcal{B}(H)$. A closed subspace $M$ of $H$ is called invariant under $T \in \mathcal{B}(H)$ if $T(M) \subseteq M$ and is said to be reducing if both $M$ and $M^\bot$ are invariant under $T$.

For a detailed study of the basic definitions and results in operator theory, we refer to \cite{CON,ReedSimon}.

\begin{definition}\label{ess_spectrum}\cite[Definiton 4.1, Page 358]{CON}
	For $T \in \mathcal{B}(H)$, the essential spectrum of $T$ is defined by
	$\sigma_{ess}(T) = \sigma(\pi(T))$, where $\pi : \mathcal{B}(H) \to \mathcal{B}(H)/ \mathcal{K}(H)$ is the canonical quotient map.
\end{definition}

In the case of a self-adjoint operator, the above definition of the essential spectrum coincides with the following description.
\begin{theorem}\label{ess spectrum of self-adjoint}\cite[Theorem VII.11, Page 236]{ReedSimon}
	Let $T = T^* \in \mathcal{B}(H)$. Then $\lambda \in \sigma_{ess}(T)$ if and only if one or more of the following conditions hold.
	\begin{enumerate}
		\item $\lambda$ is an eigenvalue of $T$ with infinite multiplicity.
		\item $\lambda$ is the limit point of $\sigma_{p}(T)$.
		\item $\lambda \in \sigma_{c}(T)$.
	\end{enumerate}
\end{theorem}

For $T \in \mathcal{B}(H)$, let $\pi_{00}(T)$ denote the set of all isolated eigenvalues of $T$ with finite multiplicity. For a self-adjoint operator $T$, $\sigma(T)\setminus \pi_{00}(T) = \sigma_{ess}(T)$.  For various equivalent definitions of the essential spectrum of a self-adjoint operator, we refer to \cite{BarrySimon4}.

For $T \in \mathcal{B}(H)$, the essential minimum modulus of $T$ \cite{BOU} is defined by
\[ m_e(T) = \inf \{\lambda : \lambda \in \sigma_{ess}(|T|)\}.\]

Let $H = H_1 \oplus H_2$ and $T \in \mathcal{B}(H)$.
Let $P_j$ be the orthogonal projection on $H$ with range $H_j$, $j=1,2$. Then $T= \begin{pmatrix}
T_{11}&T_{12}\\T_{21}&T_{22}
\end{pmatrix}$, where the operator
$ T_{i j} : H_j \to H_i$ is given by $T_{ij} =P_iTP_j|_{H_j}, i,j=1,2$.
Moreover, $T(H_1) \subseteq H_1$ if and only if $T_{12} = 0$ and $H_1$ reduces $T$ if and only if $T_{12} =0$ and $T_{21} =0$ (see \cite{CON} for more details).

Next, we define a few classes of non-normal operators.
\begin{definition} If $T \in \mathcal{B}(H)$, then $T$ is called
	\begin{enumerate}
		\item hyponormal if $\|T^{\ast} x\| \leq \|T x\|,\ x \in H.$
		\item paranormal if $\|T x\|^ 2 \leq \|T^2x\|\|x\|, \ x \in H$.
		\item $\ast$-paranormal if $\|T^*x\|^ 2 \leq \|T^2x\|\|x\|,\  x \in H$.
	\end{enumerate}
\end{definition}
For a detailed study of these classes of operators, we refer to \cite{Arora,BER1,Furutapara,HAN,Istratescu,STA,Uchi*para}.

The following result is an equivalent definition of a $\ast$-paranormal operator.
\begin{theorem} \label{*paradef}\cite{Arora,HAN} Let $T \in \mathcal{B}(H)$. Then $T$ is $\ast$-paranormal if and only if
	\[T^{*2} T^2 - 2kT T^* + k^2 I \geq 0 \ \text{for  all} \ k > 0.\]
\end{theorem}
\section{Representation and normality of $\ast$-paranormal operators}
In this section, we deduce a representation of $\ast$-paranormal operators in $\overline{\mathcal{AN}(H)}$.

Let $M := \{x \in H : \|Tx\| = \|T\|\|x\|\}$. Then by \cite[Lemma 3.1]{Rameshpara}, we have $M = N(\|T\|^2I - T^*T) = N(|T| - \|T\| I)$.

In the next result, we discuss about powers of $\ast$-paranormal operators.
\begin{theorem}\label{allpowersNA}
	If $T \in \mathcal{N}(H)$ is $\ast$-paranormal, then $T^n \in \mathcal{N}(H)$ for all $n \geq 1$.
\end{theorem}
\begin{proof}
	As $T \in \mathcal{N}(H)$, there exists a non-zero $x \in H$ such that $\|Tx\| = \|T\|\|x\|$. This is equivalent to $T^*Tx = \|T\|^2x$. Hence, we have \begin{equation}\label{inequalities} \|T\|^4\|x\|^2 = \|T^*Tx\|^2 \leq \|T^3x\|\|Tx\| \leq \|T\| \|T^2x\| \|Tx\| \leq \|T\|^2 \|Tx\|^2 = \|T\|^4 \|x\|^2.\end{equation}
	From the above inequalities, we get $\|T^2x\| = \|T\| \|Tx\| = \|T\|^2 \|x\|$ and $\|T^3x\| = \|T\| \|T^2x\| = \|T\|^3 \|x\|$.
	
	Following the similar steps as above, we conclude that $\|T^nx\| = \|T\|^n \|x\|$ for all $n \in \mathbb{N}$. As $\|T^n\| = \|T\|^n$ for all $n \in \mathbb{N}$, we get that $T^n \in \mathcal{N}(H)$ for all $n \in \mathbb{N}$.
\end{proof}

The above result need not be true if the condition of $\ast$-paranormality is dropped.

\begin{example}
	Let $T : \ell^2(\mathbb{N}) \to \ell^2(\mathbb{N})$ be defined by
	\[T(x_1,x_2,x_3,\dots) = (x_2, \left(1-\frac{1}{3}\right) x_3, \left(1-\frac{1}{4}\right) x_4, \dots), \ \forall (x_1,x_2,x_3,\dots) \in \ell^2(\mathbb{N}).\]
	
	Then
	
	\[T^2(x_1,x_2,x_3,\dots) = (\left(1-\frac{1}{3}\right)x_3, \left(1-\frac{1}{3}\right)\left(1-\frac{1}{4}\right) x_4, \dots), \ \forall (x_1,x_2,x_3,\dots) \in \ell^2(\mathbb{N}),\]
	
	\[T^*(x_1,x_2,x_3,\dots) = (0, x_1, \left(1-\frac{1}{3}\right) x_2, \left(1-\frac{1}{4}\right) x_3, \dots), \ \forall (x_1,x_2,x_3,\dots) \in \ell^2(\mathbb{N}),\]
	
	\[T^*T(x_1,x_2,x_3,\dots) = (0, x_2, \left(1-\frac{1}{3}\right)^2 x_3, \left(1-\frac{1}{4}\right)^2 x_4, \dots), \ \forall (x_1,x_2,x_3,\dots) \in \ell^2(\mathbb{N}),\]
	and \[T^{2*}T^2(x_1,x_2,x_3,\dots) = (0, 0, \left(1-\frac{1}{3}\right)^2 x_3, \left(1-\frac{1}{3}\right)^2\left(1-\frac{1}{4}\right)^2 x_4, \dots), \ \forall (x_1,x_2,x_3,\dots) \in \ell^2(\mathbb{N}).\]
	Since $1 = \|T^*(1,0,0,\dots)\|^2 \nleq 0 = \|T^2(1,0,0,\dots)\|$, we get that $T$ is not $\ast$-paranormal.
	Next, \[\|T\| = \sup\left\{1, \left(1-\frac{1}{n}\right) :n  \geq 3\right\} = 1 = \|T(0,1,0,0,\dots)\|\]
	Hence $T \in \mathcal{N}(\ell^2(\mathbb{N}))$. Also,
	\[\|T^2\| =  \sup\left\{\left(1-\frac{1}{3}\right), \left(1-\frac{1}{n}\right) \left(1-\frac{1}{n+1}\right):n  \geq 3\right\} = 1.\] Since 1 is not an eigenvalue of $T^{2*}T^2$ by \cite[Corollary 2.4]{CAR1}, we get $T^2 \notin \mathcal{N}(\ell^2(\mathbb{N}))$.
\end{example}

\begin{remark}\label{commonpointNA}
	From Theorem \ref{allpowersNA}, it is evident that if $T \in \mathcal{N}(H)$ is $\ast$-paranormal, then
	\[\{x \in H : \|T^nx\| = \|T\|^n \|x\|, n=0,1,2,\dots\} \neq \{0\}.\]
\end{remark}

In \cite{NB*-para}, the author showed the existence of an invariant subspace for a $\ast$-paranormal norm attaining operator by proving the following result.
\begin{prop}
	If $T \in \mathcal{N}(H)$ is $\ast$-paranormal, then $M_*= N(|T^*|-\|T\|I) \neq \{0\}$ and $M_ * \subseteq M$. Moreover $M_*$ is invariant under $T$.
\end{prop}

Next we show that $M$ itself is invariant under $T$.
\begin{theorem}\label{MinvunderT}
	If $T \in \mathcal{N}(H)$ is $\ast$-paranormal, then $M$ is invariant under $T$. Moreover if $M$ is finite dimensional, then $M$ reduces $T$.
\end{theorem}
\begin{proof}
	As $T \in \mathcal{N}(H)$, we have $M \neq  \{0\}$. If $x \in M$, then $\|Tx\| = \|T\| \|x\|$. From inequalities \eqref{inequalities}, we get $\|T^2x\| = \|T\| \|Tx\|$. This implies $Tx \in M$.
	
	Clearly $\frac{T|_{M}}{\|T\|}$ is an isometry. If $M$ is finite dimensional, then $T|_{M}= \|T\| U$, where $U \in \mathcal{B}(M)$ is unitary. Hence by \cite[Theorem 2.6]{RashidT*para}, $M$ reduces $T$.
\end{proof}

Now, we mainly concentrate on giving a representation for $\ast$-paranormal operators in $\overline{\mathcal{AN}(H)}$.
We first quote a few important results regarding positive operators in $\overline{\mathcal{AN}(H)}$ which will be used recurrsively in this article.

\begin{prop}\label{diagonalization}
	If $T\in \overline{\mathcal{AN}(H)}$ is positive, then $T$ is diagonalizable.
\end{prop}
\begin{theorem}\cite [Theorem 4.6]{RAMSSS}\label{positiveessentialspectrum}
	Let $T \in \mathcal{B}(H)$ be positive. Then $T \in \overline{\mathcal{AN}(H)}$  if and only if $\sigma_{ess}(T)$ is a singleton set.
\end{theorem}
Let $T \in \overline{\mathcal{AN}(H)}$ be positive and $\sigma_{ess}(T) = \{\alpha\}$, where $\alpha \geq 0$. Then  $\sigma(T) \subseteq [m(T), \alpha) \cup [\alpha, \|T\|].$ It is clear that $m_e(T)= \alpha$.

\begin{prop}\label{*para1}
	Let $T \in \overline{\mathcal{AN}(H)}$ be $\ast$-paranormal such that $\|T\|$ is an eigenvalue of $|T|$ with infinite multiplicity. Let $M = N(|T|-\|T\|I)$. Then
	\[T=\begin{blockarray}{ccc}
M & M^\bot  \\
\begin{block}{(cc)c}
\|T\| V & A & M \\0 & B & M^\perp \\
\end{block}
\end{blockarray},\]
\end{prop}
\noindent
where $V \in \mathcal{B}(M)$ is an isometry and $A \in \mathcal{B}(M^\perp, M)$, $B \in \mathcal{B}(M^\perp)$ are such that $V^*A =0$, $A^*A+B^*B \leq \|T\|^2 I_{M^\perp}$ and  $(\|T\|^4 +k^2)I_{M^\bot} -2kBB^* \geq 0$ for all $k>0$.

\begin{proof}
	Let $T = W|T|$ be the polar decomposition of $T$. Since $T \in \overline{\mathcal{AN}(H)}$, by \cite[Lemma 3.14]{RAMSSS} and Theorem \ref{positiveessentialspectrum}, we have $|T| \in \overline{\mathcal{AN}(H)}$ and $\sigma_{ess}(|T|)$ is a singleton set. As $\|T\|$ is an eigenvalue of $|T|$ with infinite multiplicity, $M \neq \{0\}$ and $\sigma_{ess}(|T|) = \{\|T\|\}$.
	
\begin{figure}[h]
\begin{tikzpicture}
\draw[very thick]  (-5,0)--(5,0);
\filldraw [gray] (5,0) circle (3pt);
\filldraw(5,0)  node[anchor=north] {$\|T\|$};
\draw [orange, very thick](-5,-.3)--(-5,.3);
\filldraw(-5,-.3) node[anchor=east] {$\beta_1$};
\draw (-.83,-.2)--(-.83,.2);
\draw (-1.12,-.2)--(-1.12,.2) node[anchor=south] {.};	
\draw (-1.45,-.2)--(-1.45,.2) node[anchor=south] {.};
\draw (-1.82,-.2)--(-1.82,.2) node[anchor=south] {.};
\draw (-2.02,-.2)--(-2.02,.2) node[anchor=south] {$\beta_5$};
\draw (-2.5,-.2)--(-2.5,.2) node[anchor=south] {$\beta_4$};
\draw (-3.12,-.2)--(-3.12,.2) node[anchor=south] {$\beta_3$};
\draw (-3.7,-.2)--(-3.7,.2) node[anchor=south] {$\beta_2$};

\draw (0,-.2)--(0,.2);
\draw (0.5,-.2)--(0.5,.2);
\draw (1,-.2)--(1,.2);
\draw (1.28,-.2)--(1.28,.2);
\draw (1.63,-.2)--(1.63,.2);
\draw (1.82,-.2)--(1.82,.2);
\draw (2.5,-.2)--(2.5,.2) ;
\draw (3.12,-.2)--(3.12,.2) ;
\draw (4.15,-.2)--(4.15,.2) ;
\draw (4.3,-.2)--(4.3,.2) ;
\draw (4.4,-.2)--(4.4,.2) ;
\draw (4.5,-.2)--(4.5,.2) ;
\draw (4.6,-.2)--(4.6,.2) ;
\draw (4.65,-.2)--(4.65,.2) ;
\draw (4.7,-.2)--(4.7,.2) ;
\draw (4.75,-.2)--(4.75,.2) ;
\draw (4.8,-.2)--(4.8,.2) ;
\draw (4.85,-.2)--(4.85,.2) ;
\draw (3.9,-.2)--(3.9,.2) ;
\end{tikzpicture}
\caption{Spectral Diagram of $|T| \in \overline{\mathcal{AN}(H)}$}
\end{figure}
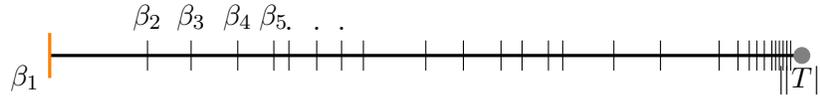

	Hence by Theorem \ref{MinvunderT}, we get $M$ is invariant under $T$. Also, for all $x \in M$, \[Wx = \frac{W|T|x}{\|T\|} = \frac{Tx}{\|T\|} \in M.\] So $M$ is also invariant under $W$. Moreover $M \subseteq N(|T|)^\bot = N(W)^\bot$.
	
Therefore, with respect to $H = M \oplus M^\perp$, $W$ and $|T|$ can be written as follows. \[W =\begin{blockarray}{ccc}
	M & M^\bot  \\
	\begin{block}{(cc)c}
	V  & V_1 & M \\0 & V_2 & M^\perp \\
	\end{block}
	\end{blockarray},
	\\
|T| = \begin{blockarray}{ccc}
	M & M^\bot  \\
	\begin{block}{(cc)c}
	\|T\| I_M  & 0& M \\0 & T_1 & M^\perp \\
	\end{block}
	\end{blockarray},\]
	where $V$ is an isometry on $M$, $V_1 = P_{M}WP_{M^\perp}|_{M^\perp}$ and  $V_2 = P_{M^\perp}WP_{M^\perp}|_{M^\perp}$ and $T_1 = \displaystyle \bigoplus_{i=1}^{n} \beta_jI_{N(|T|- \beta_jI)}$, where $\pi_{00}(|T|) = \{\beta_j\}_{i=1}^{n}$ for some $n \in \mathbb{N} \cup \{\infty\}$. If $\pi_{00}(|T|) = \emptyset$, then $|T|= \|T\|I$ and hence $T = V|T|= \|T\|V$.
	
	If $\pi_{00}(|T|) \neq \emptyset$, then
	\[T =\begin{blockarray}{ccc}
	M & M^\bot  \\
	\begin{block}{(cc)c}
	\|T\|V  & A & M \\0 & B & M^\perp \\
	\end{block}
	\end{blockarray},\]
where  $A = V_1T_1$ and $B = V_2T_1$.
Since $T^*T = |T|^2$, we get
\begin{equation*}
\begin{split}
\begin{pmatrix}
\|T\|^2 V^*V & \|T\|V^*A \\
\|T\|A^*V & A^*A +B^*B
\end{pmatrix}
=
\begin{pmatrix}
\|T\|^2 I_{M} &  0  \\
0  &  T^2_1
\end{pmatrix} \\
\end{split}.
\end{equation*}

Comparing the above two matrices,we get $V^*A =0$ and $A^*A +B^*B = T^2_1 = \displaystyle \bigoplus_{j=1}^{m}\beta^2_jI_{N(|T|-\beta_jI)} \leq \|T\|^2 I_{M^\perp}.$

Moreover,
\begin{equation*}
\hspace{-4cm}\begin{split}
TT^*
&=  \begin{pmatrix}
\|T\| V & A \\
0 & B
\end{pmatrix}
\begin{pmatrix}
\|T\| V^* & 0 \\
A^*  & B^*
\end{pmatrix} \\
&=  \begin{pmatrix}
\|T\|^2 VV^*+AA^* & AB^* \\
BA^*  & BB^*
\end{pmatrix} ,
\end{split}
\end{equation*}
and
\begin{equation*}
\begin{split}
T^{*2}T^2
&=  \begin{pmatrix}
\|T\|^2V^{*2} & 0 \\
\|T\|A^*V^*+B^*A^* & B^{*2}
\end{pmatrix}
\begin{pmatrix}
\|T\| V^2 & \|T\|VA+AB \\
0  & B^2
\end{pmatrix}\\&=  \begin{pmatrix}
\|T\|^4I_M & \|T\|^3V^*A+ \|T\|^2V^{*2}AB \\
\|T\|^3A^*V+ \|T\|^2B^*A^*V^2  & \|T\|^2A^*A+\|T\|A^*V^*AB+ \|T\|B^*A^*VA+ B^*A^*AB+B^{*2}B^2
\end{pmatrix}\\
&= \begin{pmatrix}
\|T\|^4I_M &0 \\
0& \|T\|^2A^*A+B^*A^*AB+B^{*2}B^2
\end{pmatrix}.
\end{split}
\end{equation*}

From Theorem \ref{*paradef}, we have $T$ is $\ast$-paranormal if and only if $T^{*2}T^2 -2kTT^* +k^2I \geq 0$ for all $k >0$. Substituting, we get

\begin{equation*}
\begin{split}
0 & \leq T^{*2}T^2 -2kTT^* +k^2I\\
&=  \begin{pmatrix}
(\|T\|^4+k^2)I_M -2k(\|T\|^2VV^*+AA^*) &-2kAB^* \\
-2kBA^*& \|T\|^2A^*A+B^*A^*AB+B^{*2}B^2-2kBB^*+k^2I_{M^\perp}
\end{pmatrix}
\end{split}
\end{equation*}

From the $(4,4)$ entry of the above matrix, we have \[\|T\|^2A^*A+B^*(A^*A+B^*B)B-2kBB^*+k^2I_{M^\perp} \geq 0\  \text{for all}\  k>0.\] As $A^*A +B^*B \leq \|T\|^2I_{M^\perp}$, we get \begin{equation}\label{eq1}\|T\|^4I_{M^\bot} -2kBB^*+k^2I_{M^\bot} \geq 0 \ \text{for all}\ k>0.\end{equation}\end{proof}

\begin{remark}
	\begin{enumerate}
		\item By \eqref{eq1}, it follows that $BB^* \leq \|T\|^2 I_{M^\perp}$.
		\item In Proposition \ref{*para1}, if $\pi_{00}(|T|)$ is finite, then $M^\perp$ is finite dimensional and $A \in \mathcal{B}(M^\perp,M), B \in \mathcal{B}(M^\perp,M^\perp)$ are finite rank operators.
	\end{enumerate}
\end{remark}

\begin{theorem}\label{*para2}
	Let $T \in \overline{\mathcal{AN}(H)}$ be $\ast$-paranormal. Then there exists Hilbert spaces $H_0$, $H_1$ and $H_2$ such that $H = H_0 \oplus H_1 \oplus H_2$ and $T$ can be written as \begin{equation}\label{*paraANclosure} T=\begin{blockarray}{cccc}
	H_0 & H_1& H_2  \\
	\begin{block}{(ccc)c}
V_0 & 0 & 0& H_0 \\0 & \lambda V & A&H_1 \\ 0&0&B&H_2\\
	\end{block}
	\end{blockarray},\end{equation}
where
\begin{enumerate}
	\item $\sigma_{ess}(|T|) = \{\lambda\}, \lambda \geq 0$.
	\item $H_0 = \displaystyle \bigoplus_{i=1}^{n} N(|T|-\alpha_iI)$, where $(\lambda, \|T\|] \cap \sigma(|T|) = \{\alpha_i\}_{i=1}^{n}$ for some $n \in \mathbb{N} \cup \{\infty\}$, $H_1 = N(|T|- \lambda I)$ and $H_2 = (H_0 \oplus H_1)^\perp$.
	\item $V_0 = \displaystyle \bigoplus_{i=1}^{n} \alpha_iU_i$, where $U_i \in \mathcal{B}(N(|T|-\alpha_iI))$ is unitary for all $i=1,2,\dots,n$, $n \in \mathbb{N} \cup \{\infty\}$.
	\item $V \in \mathcal{B}(H_1)$ is an isometry if $H_1 \neq \{0\}$ and $A \in \mathcal{B}(H_2,H_1 )$, $B \in \mathcal{B}(H_2)$ are such that $V^*A =0$, $A^*A+B^*B \leq \lambda^2 I_{H_2}$ and $(\lambda^4 + k^2)I_{H_2} -2kBB^* \geq 0$ for all $k>0$.
\end{enumerate}
\end{theorem}
\begin{proof}
	 To obtain the required conclusions, we apply Theorems \ref{positiveessentialspectrum} and \ref{ess spectrum of self-adjoint} to $|T|$. As $T \in \overline{\mathcal{AN}(H)}$, by \cite[Lemma 3.14]{RAMSSS}, we have $|T| \in \overline{\mathcal{AN}(H)}$ . Hence  by Theorem \ref{positiveessentialspectrum}, we have $\sigma_{ess}(|T|)$ is a singleton set, say $\{\lambda\}, \lambda \geq 0$. So $(\lambda, \|T\|] \cap \sigma(|T|)$ is at most countable say $\{\alpha_i\}_{i=1}^{n}$ for some $n \in \mathbb{N} \cup \{\infty\}$ and  $[m(T), \lambda) \cap \sigma(|T|)$ is also at most countable say $\{\beta_j\}_{j=1}^{m}$ for some $m \in \mathbb{N} \cup \{\infty\}$ . Note that $\alpha_j < \alpha_i$ and $\beta_i <\beta_j$, whenever $i <j$.
	
\begin{figure}[h]
		\begin{tikzpicture}
	\draw[orange, very thick](5,-.3)--(5,.3);
	\filldraw(5,-.3)  node[anchor=west] {$\alpha_1$};
	\draw[very thick]  (-5,0)--(5,0);
	\filldraw [gray] (0,0) circle (2.5pt);
	\filldraw(0,-.3)  node[anchor=north] {$\lambda$};
	\draw [orange, very thick](-5,-.3)--(-5,.3);
	\filldraw(-5,-.3) node[anchor=east] {$\beta_1$};
	\draw (-.11,-.2)--(-.11,.2);
	\draw (-.15,-.2)--(-.15,.2);
	\draw (-.21,-.2)--(-.21,.2);
	\draw (-.28,-.2)--(-.28,.2);
	\draw (-.37,-.2)--(-.37,.2);
	\draw (-.47,-.2)--(-.47,.2);	
	\draw (-.58,-.2)--(-.58,.2);
	\draw (-.7,-.2)--(-.7,.2);
	\draw (-.83,-.2)--(-.83,.2);
	\draw (-.97,-.2)--(-.97,.2);
	\draw (-1.12,-.2)--(-1.12,.2) node[anchor=south] {.};	
	\draw (-1.28,-.2)--(-1.28,.2) node[anchor=south] {.};
	\draw (-1.45,-.2)--(-1.45,.2) node[anchor=south] {.};
	\draw (-1.63,-.2)--(-1.63,.2) node[anchor=south] {.};
	\draw (-1.82,-.2)--(-1.82,.2) node[anchor=south] {.};
	\draw (-2.02,-.2)--(-2.02,.2) node[anchor=south] {$\beta_5$};
	\draw (-2.5,-.2)--(-2.5,.2) node[anchor=south] {$\beta_4$};
	\draw (-3.12,-.2)--(-3.12,.2) node[anchor=south] {$\beta_3$};
	\draw (-3.7,-.2)--(-3.7,.2) node[anchor=south] {$\beta_2$};
	
	\draw (.10,-.2)--(.10,.2);
	\draw (.15,-.2)--(.15,.2);
	\draw (.21,-.2)--(.21,.2);
	\draw (.28,-.2)--(.28,.2);
	\draw (.37,-.2)--(.37,.2);
	\draw (.47,-.2)--(.47,.2);	
	\draw (.58,-.2)--(.58,.2);
	\draw (.7,-.2)--(.7,.2);
	\draw (.83,-.2)--(.83,.2) node[anchor=south] {$.$};
	\draw (.97,-.2)--(.97,.2) node[anchor=south] {$.$};
	\draw (1.12,-.2)--(1.12,.2)node[anchor=south] {$.$};	
	\draw (1.28,-.2)--(1.28,.2) node[anchor=south] {$.$};
	\draw (1.45,-.2)--(1.45,.2) node[anchor=south] {$.$};
	\draw (1.63,-.2)--(1.63,.2) node[anchor=south] {$.$};
	\draw (1.82,-.2)--(1.82,.2) node[anchor=south] {$.$};
	\draw (2.02,-.2)--(2.02,.2) node[anchor=south] {$\alpha_5$};
	\draw (2.5,-.2)--(2.5,.2) node[anchor=south] {$\alpha_4$};
	\draw (3.12,-.2)--(3.12,.2) node[anchor=south] {$\alpha_3$};
	\draw (3.9,-.2)--(3.9,.2) node[anchor=south] {$\alpha_2$};
	\end{tikzpicture}
	\caption{Spectral Diagram of $|T| \in \overline{\mathcal{AN}(H)}$}
\end{figure}
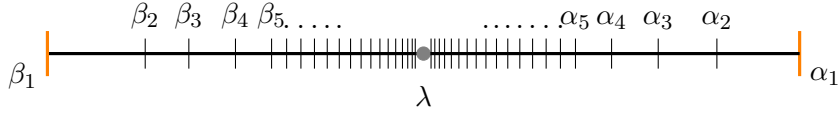

	Let $H_0 = \displaystyle \bigoplus_{i=1}^{n} N(|T|- \alpha_i I)$, $H_1 = N(|T|- \lambda I)$ and $H_2 = \displaystyle \bigoplus_{j=1}^{m} N(|T|- \beta_j I)$. If $(\lambda, \|T\|] \cap \sigma(|T|) = \emptyset$, then $H_0 = \emptyset$ and if $[m(T), \lambda) \cap \sigma(|T|) = \emptyset$, then $H_2 =\{0\}$.

Without loss of generality, let $\alpha_1 = \|T\|$. As $T$ is $\ast$-paranormal, by Theorem \ref{MinvunderT}, we have $N(|T|-\alpha_1 I)$ is invariant under $T$ and $\frac{T}{\alpha_1}$ is an isometry on $N(|T|-\alpha_1 I)$. As $\alpha_1 \in \pi_{00}(|T|)$, we have $N(|T|-\alpha_1 I)$ is finite dimensional. Hence $N(|T|-\alpha_1 I)$ reduces $T$ and $T|_{N(|T|-\alpha_1 I)} = \alpha_1U_1$, where $U_1$ is unitary on $N(|T|-\alpha_1 I)$. Next define $T_1 = T|_{N(|T|-\alpha_1 I)^\perp}$. Then $T_1$ is a $\ast$-paranormal operator and $\|T_1\| = \alpha_2$. By continuing the same process as above, we get $H_0$ reduces $T$ and $T|_{H_0} = \displaystyle \bigoplus_{i=1}^{n} \alpha_iU_i$, where $U_i$ is a unitary operator on $N(|T|-\alpha_iI)$, $i=1,2,\dots,n$, $n \in \mathbb{N} \cup \{\infty\}$.

Now, $T_2 = T|_{H_0^\perp}$ is a $\ast$-paranormal operator with $\|T_2\| = \lambda$. Then we have the following cases.

Case(1): $\lambda$ is an eigenvalue of $|T|$ of infinite multiplicity: By Proposition \ref{*para1}, $T_2$ can be written as
	\begin{equation}T_2=\begin{blockarray}{ccc}
H_1 & H_2  \\
\begin{block}{(cc)c}
\lambda V & A & H_1 \\0 & B & H_2 \\
\end{block}
\end{blockarray},\end{equation}
where $V \in \mathcal{B}(H_1)$ is an isometry and $A \in \mathcal{B}(H_2, H_1)$, $B \in \mathcal{B}(H_2)$ are such that $V^*A =0$, $A^*A+B^*B \leq \lambda^2 I_{H_2}$ and $(\lambda^4 + k^2)I_{H_2} -2kBB^* \geq 0$ for all $k>0$. Hence

\begin{equation}
T=\begin{blockarray}{cccc}
H_0 & H_1& H_2  \\
\begin{block}{(ccc)c}
V_0 & 0 & 0& H_0 \\0 & \lambda V & A&H_1 \\ 0&0&B&H_2\\
\end{block}
\end{blockarray}
\end{equation}

Subcase (a): If $\lambda$ is not a limit point of $\sigma(|T|)$, then both $H_0$, $H_2$ are finite dimensional and $V_0$, $A$ and $B$ are finite rank operators.

Subcase (b): If $\lambda$ is a limit point of $\sigma(|T|)$, then either $H_0$ or $H_2$ or both are infinite dimensional depending on $\{\alpha_i\}$ and $\{\beta_j\}$ are finite or infinite.

Case(2): $\lambda$ is an eigenvalue of $|T|$ of finite multiplicity: In this case, $\lambda$ is the limit point of $\sigma(|T|)$. Moreover, $N(|T|-\lambda I)$ reduces $T$. Hence $A=0$ and $V$ is a unitary operator on $H_1$. So $T$ can be written as

\begin{equation}
T=\begin{blockarray}{cccc}
H_0 & H_1& H_2  \\
\begin{block}{(ccc)c}
V_0 & 0 & 0& H_0 \\0 & \lambda V & 0&H_1 \\ 0&0&B&H_2\\
\end{block}
\end{blockarray}
\end{equation}

and $B$ satisfies the conditions as in Case(1). Further, Subcases (a) and (b) of Case (1) also hold here.

Case(3) : $\lambda$ is the limit point of $\sigma(T)$ but not an eigenvalue with infinite multiplicity: Then $H_1 = \{0\}$, in turn $V=0$ and $A=0$. In this case,

\begin{equation}
T=\begin{blockarray}{ccc}
H_0 & H_2  \\
\begin{block}{(cc)c}
V_0 & 0 & H_0\\ 0&B&H_2\\
\end{block}
\end{blockarray}
\end{equation}

Further, either $H_0$ or $H_2$ or both are infinite dimensional depending on $\{\alpha_i\}$ and $\{\beta_j\}$ are finite or infinite.

If $m(T) = \lambda = \|T\|$, then $T = \|T\|V$.
  \end{proof}

\begin{remark}
	If $T \in \mathcal{K}(H)$ is $\ast$-paranormal, then $T \in \overline{\mathcal{AN}(H)}$ and $\sigma_{ess}(|T|) = \{0\}$. Hence by Theorem \ref{*para2}, we get $A = B =0$ and $T$ is normal. So every $\ast$-paranormal compact operator is normal.
\end{remark}

In \cite{NB*-para}, the author has given a representation of $\ast$-paranormal $\mathcal{AN}$-operators. Here we give a representation of $\ast$-paranormal $\mathcal{AN}$-operators which is a particular case of Theorem \ref{*para2}.

\begin{theorem}\label{*paraAN1}
	If $T \in \mathcal{AN}(H)$ is  $\ast$-paranormal, then $T$ has the following representation. \begin{equation}\label{*paraAN} T=\begin{blockarray}{cccc}
	H_0 & H_1& H_2  \\
	\begin{block}{(ccc)c}
	V_0 & 0 & 0& H_0 \\0 & \lambda V & A&H_1 \\ 0&0&B&H_2\\
	\end{block}
	\end{blockarray},\end{equation}
	where
	\begin{enumerate}
		\item $\sigma_{ess}(|T|) = \{\lambda\}, \lambda \geq 0$.
		\item $H_0 = \displaystyle \bigoplus_{i=1}^{n} N(|T|-\alpha_iI)$, where $(\lambda, \|T\|] \cap \sigma(|T|) = \{\alpha_i\}_{i=1}^{n}$ for some $n \in \mathbb{N} \cup \{\infty\}$, $H_1 = N(|T|- \lambda I)$ and $H_2 = (H_0 \oplus H_1)^\perp$ is finite dimensional.
		\item $V_0 = \displaystyle \bigoplus_{i=1}^{n} \alpha_iU_i$, where $U_i \in \mathcal{B}(N(|T|-\alpha_iI))$ is unitary for all $i=1,2,\dots,n$, $n \in \mathbb{N} \cup \{\infty\}$.
		\item $V \in \mathcal{B}(H_1)$ is an isometry if $H_1 \neq \{0\}$ and $A \in \mathcal{B}(H_2,H_1 )$, $B \in \mathcal{B}(H_2)$ are finite rank operators such that $V^*A =0$, $A^*A+B^*B \leq \lambda^2 I_{H_2}$ and $(\lambda^4 + k^2)I_{H_2} -2kBB^* \geq 0$ for all $k>0$.
	\end{enumerate}
\end{theorem}
\begin{proof}
	By \cite[Corollaries 2.10,2.11]{VEN}, we have $T \in \mathcal{AN}(H)$ if and only if $|T| \in \mathcal{AN}(H)$. Now by \cite[Theorem 2.4]{Rameshpara}, we have $|T| \in \mathcal{AN}(H)$ if and only if $\sigma_{ess}(|T|)$ is a singleton set and $[m(T),m_e(T)) \cap \sigma(|T|)$ contains only finitely many eigenvalues of $|T|$. So by the similar steps as in Theorem \ref{*para2}, we get the required representation.
\end{proof}

\begin{remark}
  We note that in Theorem \ref{*paraAN1}, the space $H_2$ is finite dimensional, which guarantees that the operators $A$ and $B$ have finite rank, which is not the case in \cite{NB*-para}.
\end{remark}

Since $\mathcal{AM}$-operators are contained in  $\overline{\mathcal{AN}(H)}$, using Theorem \ref{*para2}, the following representation of $\mathcal{AM}$-operators is obtained.

\begin{theorem}\label{*paraAM1}
	If $T \in \mathcal{AM}(H)$ is  $\ast$-paranormal, then $T$ has the following representation. \begin{equation}\label{*paraAM} T=\begin{blockarray}{cccc}
		H_0 & H_1& H_2  \\
		\begin{block}{(ccc)c}
		V_0 & 0 & 0& H_0 \\0 & \lambda V & A&H_1 \\ 0&0&B&H_2\\
		\end{block}
		\end{blockarray},\end{equation}
		where
		\begin{enumerate}
			\item $\sigma_{ess}(|T|) = \{\lambda\}, \lambda \geq 0$.
			\item $H_0 = \displaystyle \bigoplus_{i=1}^{n} N(|T|-\alpha_iI)$, where $(\lambda, \|T\|] \cap \sigma(|T|) = \{\alpha_i\}_{i=1}^{n}$ for some $n \in \mathbb{N}$ is finite dimensional, $H_1 = N(|T|- \lambda I)$ and $H_2 = (H_0 \oplus H_1)^\perp$.
			\item $V_0 = \displaystyle \bigoplus_{i=1}^{n} \alpha_iU_i$ is a finite rank normal operator, where $U_i \in \mathcal{B}(N(|T|-\alpha_iI))$ is unitary for all $i=1,2,\dots,n$, $n \in \mathbb{N} \cup \{\infty\}$.
			\item $V \in \mathcal{B}(H_1)$ is an isometry if $H_1 \neq \{0\}$ and $A \in \mathcal{B}(H_2,H_1 )$, $B \in \mathcal{B}(H_2)$ are such that $V^*A =0$, $A^*A+B^*B \leq \lambda^2 I_{H_2}$ and $(\lambda^4 + k^2)I_{H_2} -2kBB^* \geq 0$ for all $k>0$.
		\end{enumerate}
	\end{theorem}
\begin{proof}
	By \cite[Theorem 5.14]{GAN}, we have $T \in \mathcal{AM}(H)$ if and only if $|T| \in \mathcal{AM}(H)$. Now by \cite[Theorem 3.10]{BALA1}, $|T| \in \mathcal{AM}(H)$ if and only if $\sigma_{ess}(|T|)$ is a singleton set and $(m_e(T), \|T\|] \cap \sigma(|T|)$ contains only finitely many eigenvalues of $|T|$. So by the similar steps as in Theorem \ref{*para2}, we get the required representation.
\end{proof}

Next, we want to give a representation of $T \in \overline{\mathcal{AN}(H)}$ when $T^*$ is $\ast$-paranormal.

\begin{corollary}\label{T^*isANclosure}
	If $T \in \overline{\mathcal{AN}(H)}$ is such that $T^*$ is $\ast$-paranormal, then $T^* \in \overline{\mathcal{AN}(H)}$.
\end{corollary}

\begin{proof}
	As $T \in \overline{\mathcal{AN}(H)}$, by \cite[Theorem 3.15]{RAMSSS}, we have  $T^*T \in \overline{\mathcal{AN}(H)}$ and by Theorem \ref{positiveessentialspectrum}, $\sigma_{ess}(T^*T) = \{\lambda\}$ for some $\lambda \geq 0$.
	
	By \cite[Theorem 6, Page 173]{MUL}, it is well known that $\sigma_{ess}(T^*T) \setminus \{0\} = \sigma_{ess}(TT^*)\setminus \{0\}$.
	
	If $\lambda =0$, then $T^*T \in \mathcal{K}(H)$, so is $TT^*$. Hence in this case, $\sigma_{ess}(T^*T) = \sigma_{ess}(TT^*)$.
	
	So let us consider $\lambda \neq 0$. If $0 \in \sigma_{ess}(TT^*)$, then by \cite[Theorem VII.11, Page 236]{ReedSimon} either  $0$ is an eigenvalue of $TT^*$ of infinite multiplicity or it is the limit point of $\sigma(TT^*)$.
	
	If $0$ is an eigenvalue of $TT^*$ of infinite multiplicity, then $N(T^*)$ is infinite dimensional as $N(TT^*) = N(T^*)$. By the $\ast$-paranormality of $T^*$, we have $N(T^*) \subseteq N(T^{*2})\subseteq N(T) = N(T^*T)$. Hence $0 \in \sigma_{ess}(T^*T)$, which is a contradiction to $\lambda \neq 0$.
	
	If $0$ is a limit point of $\sigma(TT^*)$, then there exists a sequence $\{\lambda_n\} \subseteq \sigma(TT^*) \setminus \{0\}$ such that $\lambda_n \to 0$ as $n \to \infty$.  Since $\sigma(TT^*) \setminus \{0\} = \sigma(T^*T) \setminus \{0\}$, we get $0 \in \sigma_{ess}(T^*T)$ contradicting our assumption. Therefore $\sigma_{ess}(T^*T) = \sigma_{ess}(TT^*)$. Hence by \cite[Theorem 4.13]{RAMSSS1}, it follows that $T^* \in \overline{\mathcal{AN}(H)}$. \end{proof}

\begin{corollary}
	
Let $T \in \overline{\mathcal{AN}(H)}$ be such that $T^*$ is $\ast$-paranormal. Then there exists Hilbert spaces $H_0, H_1,H_2$ such that $H = H_0 \oplus H_1 \oplus H_2$ and $T$ can be written as \begin{equation} T=\begin{blockarray}{cccc}
	H_0 & H_1& H_2  \\
	\begin{block}{(ccc)c}
	V_0 & 0 & 0& H_0 \\0 & \lambda V & 0&H_1 \\ 0&A&B&H_2\\
	\end{block}
	\end{blockarray},\end{equation}
	where
	\begin{enumerate}
		\item $\sigma_{ess}(|T^*|) = \{\lambda\}, \lambda \geq 0$.
		\item $H_0 = \displaystyle \bigoplus_{i=1}^{n} N(|T^*|-\alpha_iI)$, where $(\lambda, \|T\|] \cap \sigma(|T^*|) = \{\alpha_i\}_{i=1}^{n}$ for some $n \in \mathbb{N} \cup \{\infty\}$, $H_1 = N(|T^*|- \lambda I)$ and $H_2 = (H_0 \oplus H_1)^\perp$.
		\item $V_0 = \displaystyle \bigoplus_{i=1}^{n} \alpha_iU_i$, where $U_i \in \mathcal{B}(N(|T^*|-\alpha_iI))$ is unitary for all $i=1,2,\dots,n$, $n \in \mathbb{N} \cup \{\infty\}$.
		\item $V \in \mathcal{B}(H_1)$ is a co-isometry if $H_2 \neq \{0\}$ and $A \in \mathcal{B}(H_1,H_2 )$, $B \in \mathcal{B}(H_2)$ are such that $VA^* =0$, $AA^*+BB^* \leq \lambda^2 I_{H_2}$ and $(\lambda^4 + k^2)I_{H_2} -2kB^*B \geq 0$ for all $k>0$.
	\end{enumerate}
\end{corollary}
\begin{proof}
	The proof follows by Corollary \ref{T^*isANclosure} and Theorem \ref{*para2}.
\end{proof}

\begin{prop}\label{4thentrypositive}
	If $T = \begin{pmatrix}
	C &A\\0&B\end{pmatrix}$ is hyponormal, then $A^*A+B^*B-BB^* \geq 0$ and $C^*C \geq CC^* +AA^*$. In particular, $C$ is hyponormal.
\end{prop}
\begin{proof}
	By the definition of hyponormality of $T$, we have
	\begin{equation*}
	\begin{split}
	0 & \leq  T^*T -TT^*\\
	&=  \begin{pmatrix}
	C^*C- CC^*-AA^* & C^*A-AB^* \\
	A^*C-BA^*  & A^*A +B^*B-BB^*
	\end{pmatrix}
	\end{split}
	\end{equation*}
Since the (1,1) and (4,4) entries of the matrix are positive, hence we get the required result.
\end{proof}
In \cite{RAMSSS1}, a similar kind of representation is also given for hyponormal operators in $\overline{\mathcal{AN}(H)}$. This need not imply every $\ast$-paranormal operator in $\overline{\mathcal{AN}(H)}$ is hyponormal.

\begin{example}
	Let $T: \ell^2(\mathbb{N}) \to \ell^2(\mathbb{N})$ be defined by
	\[T(x_1,x_2,x_3,\dots) = (0, x_1,\sqrt{2}x_2, 2x_3,2x_4,\dots), \ \forall (x_n) \in \ell^2(\mathbb{N})\]
and $S: \ell^2(\mathbb{N}) \to \ell^2(\mathbb{N})$ be defined by
\[S(x_1,x_2,x_3,\dots) = (0, \sqrt{2}x_1,x_2, 2x_3,2x_4,\dots), \ \forall (x_n) \in \ell^2(\mathbb{N}).\]
Then \[T^*(x_1,x_2,x_3,\dots) = (x_2,\sqrt{2}x_3, 2x_4,2x_5,\dots), \ \forall (x_n) \in \ell^2(\mathbb{N})\] and

\[S^*(x_1,x_2,x_3,\dots) = (\sqrt{2}x_2,x_3, 2x_4,2x_5,\dots), \ \forall (x_n) \in \ell^2(\mathbb{N}).\]
Also, \[T^*T(x_1,x_2,x_3,\dots) = (x_1,2x_2, 4x_3,4x_4,\dots), \ \forall (x_n) \in \ell^2(\mathbb{N}),\]
\[S^*S(x_1,x_2,x_3,\dots) = (2x_1,x_2, 4x_3,4x_4,\dots), \ \forall (x_n) \in \ell^2(\mathbb{N}).\]

Clearly $\|T^*x\| \leq \|Tx\|$ for all $x \in \ell^2(\mathbb{N})$, but $\|S^*(0,1,0,\dots)\| \geq \|S(0,1,0,\dots)\|$. Hence $T$ is hyponormal and $S$ is not hyponormal, but $S$ is $\ast$-paranormal by \cite[Example 1]{Uchi*para}. Also, $\sigma_{ess}(T^*T) = \sigma_{ess}(S^*S) =\{4\}$. Hence by Theorem \ref{positiveessentialspectrum}, $T^*T, S^*S \in \overline{\mathcal{AN}(H)}$ which in turn implies $T, S \in \overline{\mathcal{AN}(H)}$ \cite[Theorem 3.15]{RAMSSS}.

Let $\{e_n\}$ be the standard orthonormal basis of $\ell^2(\mathbb{N})$. Clearly $M = N(|T|-\|T\|I) = N(S-\|S\|I) = \overline{\spn}\{e_3,e_4,\dots\}$. Then
\[T=\begin{blockarray}{ccc}
M & M^\bot  \\
\begin{block}{(cc)c}
2 R & A_1 & M \\0 & B_1 & M^\perp \\
\end{block}
\end{blockarray}
\\
\\
\hspace{1cm}S=\begin{blockarray}{ccc}
M & M^\bot  \\
\begin{block}{(cc)c}
2R & A_2 & M \\0 & B_2 & M^\perp \\
\end{block}
\end{blockarray}, \]

where $R(x_1,x_2,\dots) = (0,x_1,x_2,\dots)$ for all $(x_n) \in M$,

 $A_1(x_1,x_2) = (\sqrt{2}x_2,0,\dots)$ for all $(x_1,x_2) \in M^\perp$, $B_1(x_1,x_2) = (0,x_1)$ for all $(x_1,x_2) \in M^\perp.$

 $A_2(x_1,x_2) = (x_2,0,\dots)$ for all $(x_1,x_2) \in M^\perp$, $B_2(x_1,x_2) = (0,\sqrt{2}x_1)$ for all $(x_1,x_2) \in M^\perp.$

 By Proposition \ref{4thentrypositive}, we have $\|A_1x\|^2 \geq \|B^*_1x\|^2 -\|B_1x\|^2$ for all $x \in M^\perp$. But $1= \|A_2e_2\|^2 < \|B^*_2e_2\|^2 -\|B_2e_2\|^2= 2$.
\end{example}

Next we look at some sufficient conditions implying the normality of $\ast$-paranormal operators in $\overline{\mathcal{AN}(H)}$.	

\begin{theorem}\label{inv*paraANclosure}
Let $T \in \overline{\mathcal{AN}(H)}$ be a $\ast$-paranormal operator. If $T$ is invertible, then $T$ is normal.
\end{theorem}

\begin{proof}
Let $X$ denote the set of all $\ast$-paranormal operators in $\mathcal{B}(H)$. Then $X = \cap_{k > 0}\{T : T^{*2}T^2- 2kTT^*+k^2 I \geq 0\}$ is a closed set. Hence the set of all $\ast$-paranormal operators in $\overline{\mathcal{AN}(H)}$ is closed. As the set of all invertible operators in $\mathcal{B}(H)$ is open, if $T \in \overline{\mathcal{AN}(H)}$ is an invertible, $\ast$-paranormal operator, there exists a sequence of invertible $\ast$-paranormal $\mathcal{AN}$-operators $\{T_n\}$ such that $T_n \to T$ in the operator norm. By \cite[Theorem 3.14]{NB*-para}, we have $T_n$ is normal for all $n \in \mathbb{N}$. Hence $T$ is normal.
\end{proof}

We further improve Theorem \ref{inv*paraANclosure} as follows:

\begin{corollary}\label{N(T)=N(T^*)Tisnormal}
	Let $T \in \overline{\mathcal{AN}(H)}$ be a $\ast$-paranormal operator. If $N(T) = N(T^*)$, then $T$ is normal.
	\end{corollary}
\begin{proof}
	If $T \in \overline{\mathcal{AN}(H)}$ is compact, then by \cite[Theorem 4.6]{Rashidn*}, $T$ is normal. If $T$ is not compact, then by \cite[Corollary 3.8]{RAMSSS}, $N(T)$ is finite dimensional and $R(T)$ is closed. As $N(T) = N(T^*)$, we have $N(T)$ reduces $T$. Hence
	\[T =\begin{blockarray}{ccc}
	N(T) & N(T)^{\bot}  \\
	\begin{block}{(cc)c}
	0  & 0 & N(T) \\0 & T_1 & N(T)^\perp \\
	\end{block}
	\end{blockarray},\]
	where $T_1 \in \overline{\mathcal{AN}(N(T)^\bot)}$ is an invertible, $\ast$-paranormal operator. Hence by Theorem \ref{inv*paraANclosure}, $T_1$ is normal which in turn implies that $T$ is normal.
\end{proof}

\begin{corollary}
	Let $T \in \overline{\mathcal{AN}(H)}$ be a $\ast$-paranormal operator. If dim $N(T)$ = dim $N(T^*)$, then $T$ is normal.
\end{corollary}
\begin{proof}
	If $T$ is $\ast$-paranormal, then $N(T) \subseteq N(T^2) \subseteq N(T^*)$. So dim $N(T)$ = dim $N(T^*)$ implies that $N(T) = N(T^*)$. Hence by Corollary \ref{N(T)=N(T^*)Tisnormal}, $T$ is normal.
\end{proof}

\begin{corollary}
	\begin{enumerate}
		\item If $T \in \mathcal{AN}(H)$ be a $\ast$-paranormal operator with dim $N(T)$ = dim $N(T^*)$, then $T$ is normal.
		\item If $T \in \mathcal{AM}(H)$ be a $\ast$-paranormal operator with dim $N(T)$ = dim $N(T^*)$, then $T$ is normal.
	\end{enumerate}
\end{corollary}

\section{$\ast$-paranormal Toeplitz and Hankel operators in $\overline{\mathcal{AN}(H^2)}$}
In this section, we discuss a few results on concrete operators like Toeplitz and Hankel operators because of their application in many fields. Toeplitz and Hankel operators in  $\overline{\mathcal{AN}(H^2)}$ are described in \cite{RAMSSSJMAA,RAMSSSAMToep}. In this context, we investigate such operators with the additional assumption of $ast$-paranormality.

 Let $L^2$ denote the space of all square integrable functions on the unit circle $\mathbb{T}$ in the complex plane with respect to the normalized Lebesgue measure $\mu$. For any essentially bounded measurable function $\varphi \in L^\infty$ on the circle, we have the Laurent operator $L_\varphi : L^2 \to L^2$ defined by \[L_\varphi f(z) = \varphi(z)f(z),\  \forall f \in L^2, z \in \mathbb{T}.\]
 Let $H^2$ denote the closed subspace of $L^2$ consisting of all those functions with vanishing negative Fourier co-efficients. The Laurent operator gives rise to two more operators called the Toeplitz operator $T_\varphi :H^2 \to H^2$ defined by
\[T_\varphi f = P L_\varphi f \quad  \text{for all} \ f  \in H^2,\]
and the Hankel operator  $H_\varphi :H^2 \to H^2$ defined by
\[H_\varphi f = J(I-P)L_\varphi f \quad \text{for all} \ f  \in H^2,\]
where $P$ is the orthogonal projection of $L^2$ onto $H^2$ and $J$ is the unitary operator on $L^2$ defined by $J(z^{-n}) = z^{n-1}, n= 0,\pm1,\pm2, \dots$. It is well known that $\|L_\varphi\| = \|T_\varphi\| = \|\varphi\|_\infty$.

 We prove that a $\ast$-paranormal Toeplitz operator in $\overline{\mathcal{AN}(H^2)}$ is a scalar multiple of an isometry and a $\ast$-paranormal Hankel operator in $\overline{\mathcal{AN}(H^2)}$ has to be normal.

\begin{theorem} \label{*paraNAToep}
Let $\varphi \in L^\infty$. If $T_\varphi \in \mathcal{N}(H^2)$ is  $\ast$-paranormal, then $T_\varphi$ is a scalar multiple of an isometry.
\end{theorem}
\begin{proof}
From Remark \ref{commonpointNA}, we have \[\{f \in H^2: \|T^n_\varphi f\| = \|T_\varphi\|^n\|f\|, n=0,1,2,\dots,\} \neq \{0\}.\] Hence by \cite[Theorem 6]{YOSnormatt}, $T_\varphi$ is a scalar multiple of an isometry .
\end{proof}

\begin{remark}
Let $\varphi \in L^\infty$ be such that $T_\varphi \in \mathcal{N}(H^2)$. Then $T_\varphi$ is $\ast$-paranormal if and only if $T_\varphi$ is a scalar multiple of an isometry.
\end{remark}

By \cite{RAMSSS}, it is well-known that the class of norm attaining operators neither contains the closure of $\mathcal{AN}$-operators nor contained in it. Hence we characterize all $\ast$-paranormal Toeplitz operators in $\overline{\mathcal{AN}(H^2)}$.

\begin{corollary}\label{toeplANclo}
	Let $\varphi \in L^\infty$. If $T_\varphi \in \overline{\mathcal{AN}(H^2)}$ is $\ast$-paranormal, then $T_\varphi$ is a scalar multiple of an isometry.
\end{corollary}
\begin{proof}
	Let $T_\varphi \in \overline{\mathcal{AN}(H^2)}$. Since $\|T_\varphi\| = \|\varphi\|_\infty$ for all $\varphi \in L^\infty$, there exists a sequence of $\ast$-paranormal $\mathcal{AN}$-Toeplitz operators $\{T_{\psi_n}\}$ such that $T_{\psi_n}$ converges to $T_\varphi$ in the operator norm. Since $\{T_{\psi_n}\} \subseteq \mathcal{N}(H^2)$, by Theorem \ref{*paraNAToep}, we have $T_{\psi_n} = \alpha_n T_{\varphi_n}$, where $T_{\varphi_n}$ is an isometry for all $n \in \mathbb{N}$ and $\alpha_n \in \mathbb{C}$. As $T^*_{\psi_n}T_{\psi_n} = |\alpha_n|^2$ converges to $T^*_\varphi T_\varphi$, we get $\{\alpha_n\}_{n \in \mathbb{N}}$ is bounded and hence it has a convergent subsequence which converges to some $\alpha \in \mathbb{C}$. Hence $T^*_\varphi T_\varphi = |\alpha|^2$. This implies that $T_\varphi$ is a scalar multiple of an isometry.
\end{proof}

\begin{remark}
	Let $\varphi \in L^\infty$ be such that $T_\varphi$ is $\ast$-paranormal. Then $T_\varphi \in \mathcal{N}(H^2)$ if and only if $\overline{\mathcal{AN}(H^2)}$.
\end{remark}
\begin{proof}
	The proof follows from Theorem \ref{*paraNAToep} and Corollary \ref{toeplANclo}.
\end{proof}

For Hankel operators we have the following result.

\begin{theorem}\label{*paraNAHank}
Let $\varphi \in L^\infty$.	If $H_\varphi \in \mathcal{N}(H^2)$ is $\ast$-paranormal, then $H_\varphi$ is normal.
\end{theorem}
\begin{proof}
By Remark \ref{commonpointNA} and \cite[Theorem 7]{YOSnormatt}, it follows that $H_\varphi$ is normal.
\end{proof}

\begin{remark}
	Let $\varphi \in L^\infty$ be such that $H_\varphi \in \mathcal{N}(H^2)$. Then $H_\varphi$ is $\ast$-paranormal if and only if $H_\varphi$ is normal.
\end{remark}

\begin{corollary}\label{HankelANcl}
	Let $\varphi \in L^\infty$. If $H_\varphi \in \overline{\mathcal{AN}(H^2)}$ is $\ast$-paranormal, then $H_\varphi$ is normal.
\end{corollary}
\begin{proof}
	Let $H_\varphi \in \overline{\mathcal{AN}(H^2)}$. As $\|H_\varphi\| \leq \|\varphi\|_\infty$ for all $\varphi \in L^\infty$, there exists a sequence of $\ast$-paranormal $\mathcal{AN}$-Hankel operators $\{H_{\varphi_n}\}$ such that $\{H_{\varphi_n}\}$ converges to $H_\varphi$ in the operator norm. As $\{H_{\varphi_n}\} \subseteq \mathcal{N}(H^2)$, by Theorem \ref{*paraNAHank}, we have $H_{\varphi_n}$ is normal for all $n \in \mathbb{N}$. Hence $H_\varphi$ is normal.
\end{proof}

\begin{remark}
	Let $\varphi \in L^\infty$ be such that $H_\varphi$ is $\ast$-paranormal. Then $H_\varphi \in \mathcal{N}(H^2)$ if and only if $H_\varphi \in \overline{\mathcal{AN}(H^2)}$.
\end{remark}
\begin{proof}
	The proof follows from Theorem \ref{*paraNAHank} and Corollary \ref{HankelANcl}.
\end{proof}
\section*{Funding}
The first author is supported by SERB Grant No. MTR/2019/001307, Govt. of India. The second author is supported by the  Department of Science and Technology- INSPIRE Fellowship (Grant No. DST/INSPIRE FELLOWSHIP/2018/IF180107).

\section*{Data Availability}
Data sharing is not applicable to this article as no datasets were generated or
analyzed during the current study.

\section*{Conflict of interest}
The authors declare that there are no conflicts of interests.

\end{document}